\documentclass[11pt]{amsart}

\usepackage{amsmath,amsthm}

\newtheorem{theorem}{Theorem}[section]
\newtheorem{lemma}[theorem]{Lemma}
\newtheorem{corollary}[theorem]{Corollary}
\newtheorem{prop}[theorem]{Proposition}

\theoremstyle{definition}
\newtheorem{definition}[theorem]{Definition}

\newtheorem{ass}[theorem]{Assumption}

\theoremstyle{remark}
\newtheorem{remark}[theorem]{Remark}

\numberwithin{equation}{section}

\def\R{\mathbb{R}}
\def\co{C_0(\R^n)}
\def\cb{C_b(\R^n)}
\def\e{\rm e}
\def\u{\rm u}

\def\N{\mathbb{N}}
\def\P{\mathbb{P}}

\def\rp{\mathbf{p}}
\def\y{\mathbf{y}}
\def\rq{\mathbf{q}}

\def\mP{\mathcal{P}}
\def\K{\mathbb{K}}
\def\Q{\mathbb{Q}}
\def\C{\mathbb{C}}
\def\X{\mathcal{X}}
\def\Y{\mathcal{Y}}
\def\Z{\mathcal{Z}}
\def\W{\mathcal{W}}

\def\b{\mathcal{B}}
\def\ring{\R[x]}
\def\ip#1#2{\langle#1,#2\rangle}

\begin{document}

\title[sos approximation of nonnegative polynomials]{S.o.s. approximation of polynomials,
nonnegative on a real algebraic set}

\author{Jean B. Lasserre}

\address{LAAS-CNRS\\
7 Avenue du Colonel Roche\\
31077 Toulouse C\'{e}dex 4, France}

\subjclass{11E25 12D15 13P05, 12Y05, 90C22, 90C25}

\date{}

\keywords{Global optimization; real algebraic geometry; semidefinite relaxations}

\begin{abstract}
Wih every real polynomial $f$, we associate a family 
$\{f_{\epsilon r}\}_{\epsilon, r}$ of real polynomials, 
in {\it explicit} form in terms of $f$ and the parameters $\epsilon>0,r\in\N$, 
and such that $\Vert f-f_{\epsilon r}\Vert_1\to 0$ as
$\epsilon\to 0$.

Let $V\subset \R^n$ be a real algebraic set
described by finitely many polynomials equations $g_j(x)=0,j\in J$,
and let $f$ be a real polynomial, nonnegative on $V$.
We show that for every $\epsilon>0$, there exist
nonnegative scalars $\{\lambda_j(\epsilon)\}_{j\in J}$
such that, for all $r$ sufficiently large,
\[f_{\epsilon r}+\sum_{j\in J} \lambda_j(\epsilon)\,g_j^2,\quad\mbox{is
a sum of squares.}\]
This representation is an obvious certificate of nonnegativity 
of $f_{\epsilon r}$ on $V$, and very specific 
in terms of the $g_j$ that define the set $V$. In particular, it
is valid with {\it no} assumption on $V$.
In addition, this representation is also useful from a computation point of view, 
as we can define semidefinite programing relaxations to approximate
the global minimum of $f$ on a real algebraic set $V$, or a
semi-algebraic set $\K$, and again, with {\it no} assumption on $V$ or $\K$.
\end{abstract}

\maketitle

\section{Introduction}

Let $V\subset\R^n$ be the real algebraic set
\begin{equation}
\label{setv}
V\,:=\,\{x\in\R^n\,\vert\, \quad g_j(x)\,=\,0,\quad j=1,\ldots,m\},
\end{equation}
for some family of real polynomials $\{g_j\}\subset\ring (=\R[x_1,\ldots,x_n])$.

The main motivation of this paper is to provide a characterization of polynomials
$f\in\ring$, nonnegative on $V$, in terms of a {\it certificate} of
positivity. In addition, and in view of the
many potential applications, one would like to obtain a representation
that is also {\it useful} from a computational point of view.

In some particular cases, when $V$ is {\it compact}, and viewing the
equations $g_j(x)=0$ as two opposite inequations $g_j(x)\geq0$ and $g_j(x)\leq0$,
one may obtain Schm\"udgen's  {\it sum of squares} (s.o.s.)
representation \cite{schmudgen}
for $f+\epsilon$ ($\epsilon >0$), instead of $f$. Under an additional assumption on
the $g_j$'s that define $V$,
the latter representation may be even refined to become
Putinar \cite{putinar} and Jacobi and Prestel \cite{jacobi}
s.o.s. representation, that is, $f+\epsilon$ can be written
\begin{equation}
\label{putinar}
f+\epsilon\,=\,f_0+ \sum_{j=1}^mf_j\, g_j,
\end{equation}
for some polynomials  $\{f_j\}\subset\ring$, with $f_0$ a s.o.s.
Hence, if $f$ is nonnegative on $V$, every approximation $f+\epsilon$
of $f$ (with $\epsilon >0$) has the representation (\ref{putinar}).
The interested reader is referred to 
Marshall \cite{marshall}, Prestel and Delzell \cite{prestel}, and Scheiderer \cite{claude,claude2} for a nice account of such results.

{\bf Contribution.} We propose the following result: 
Let $\Vert f\Vert_1=\sum_\alpha \vert f_\alpha\vert$ 
whenever $x\mapsto f(x)=\sum_\alpha f_\alpha x^\alpha$).
Let $f\in\ring$ be nonnegative on $V$, as defined in (\ref{setv}), and
let $F:=\{f_{\epsilon r}\}_{\epsilon, r}$ be the family of polynomials
\begin{equation}
\label{fepsilon}
f_{\epsilon r}\,=\,f+\epsilon\sum_{k=0}^{r}\sum_{i=1}^n\frac{x^{2k}_i}{k{\rm
!}},\qquad \epsilon\geq 0,\quad r\in\N.
\end{equation}
(So, for every 
$r\in\N$, $\Vert f-f_{\epsilon r}\Vert_1\to 0$ as $\epsilon\downarrow 0$.)

Then, for every $\epsilon>0$, there exist
nonnegative {\it scalars} $\{\lambda_j(\epsilon)\}_{j=1}^m$, such that
for all $r$ sufficiently large (say $r\geq r(\epsilon)$),
\begin{equation}
\label{intro-1}
f_{\epsilon r}\,=\,q_\epsilon-\sum_{j=1}^m\,\lambda_j(\epsilon)\,g_j^2,
\end{equation}
for some s.o.s. polynomial $q_\epsilon\in\ring$, that is,
$f_{\epsilon r}+\sum_{j=1}^m\lambda_j(\epsilon)g_j^2$ is s.o.s.

Thus, with {\it no} assumption on the set $V$, one obtains a
representation of $f_{\epsilon r}$ (which is positive on $V$ as 
$f_{\epsilon r}>f$ for all $\epsilon>0$) in the simple and explicit form
(\ref{intro-1}), an obvious {\it certificate} of positivity
of $f_{\epsilon r}$ on $V$.
In particular, when $V\equiv \R^n$, one retrieves
the result of \cite{lasserre}, which states that every nonnegative
real polynomial $f$ can be aproximated as closely as desired, by a
family of s.o.s. polynomials $\{f_{\epsilon r(\epsilon)}\}_\epsilon$,
with $f_{\epsilon r}$ as in (\ref{fepsilon}).

Notice that $f+n\epsilon =f_{\epsilon 0}$. 
So, on the one hand, the approximation $f_{\epsilon r}$ in (\ref{intro-1}) is more
complicated than $f+\epsilon$ in (\ref{putinar}),
valid for the compact case with an additional assumption, 
but on the other hand,
the coefficients of the $g_j$'s in (\ref{intro-1}) are now {\it scalars} instead of s.o.s., and
(\ref{intro-1}) is valid for an arbitrary algebraic set $V$.

The case of a semi-algebraic set $\K=\{x\in\R^n\vert g_j(x)\geq 0,\:j=1,\ldots ,m\}$ 
reduces to the case of an algebraic set $V\in\R^{n+m}$, 
by introducing $m$ {\it slack} variables $\{z_j\}$, 
and replacing $g_j(x)\geq 0$ with $g_j(x)-z_j^2=0$, for all $j=1,\ldots,m$.
Let $f\in\ring$ be nonnegative on $\K$. Then, for every $\epsilon >0$,
there exist nonnegative scalars $\{\lambda_j(\epsilon)\}_{j=1}^m$
such that, for all sufficiently large $r$,
\[f+\epsilon\sum_{k=0}^{r}\left[\sum_{i=1}^n\frac{x^{2k}_i}{k{\rm!}}
+\sum_{j=1}^m\frac{z^{2k}_j}{k{\rm!}}\right]\,=\,
q_\epsilon-\sum_{j=1}^m\,\lambda_j(\epsilon)\,(g_j-z_j^2)^2,\] 
for some s.o.s. $q_\epsilon\in\R[x,z]$.
Equivalently, everywhere on $\K$, the polynomial
\[x\mapsto f(x)+\epsilon\sum_{k=0}^{r}\sum_{i=1}^n\frac{x^{2k}_i}{k{\rm
!}},+\epsilon\sum_{k=0}^{r}\sum_{j=1}^m\frac{g_j(x)^k}{k{\rm !}}\]
coincides with the polynomial $x\mapsto q_\epsilon(x_1,\ldots,x_n,\sqrt{g_1(x)},\ldots,\sqrt{g_m(x)})$, obviously nonnegative.

The representation (\ref{intro-1}) is also useful for computational
purposes. Indeed, 
using (\ref{intro-1}), one can approximate the global minimum of $f$ on $V$, by
solving a sequence of semidefinite programming (SDP) problems. The same
applies to an arbitrary semi-algebraic set $\K\subset\R^n$, defined by $m$
polynomials inequalities, as explained above.
Again, and in contrast to previous SDP-relaxation techniques as in
e.g. \cite{lasserre1,lasserre2,lasserre3,parrilo2,markus}, no
compacity assumption on $V$ or $\K$ is required.

In a sense, the family $F=\{f_{\epsilon r}\}\subset\ring$ (with $f_{0
r}\equiv f$) is a set of
{\it regularizations} of $f$, because one may approximate $f$ by 
members of $F$, and those members {\it always} have nice
representations when $f$ is nonnegative on an algebraic set $V$ 
(including the case $V\equiv \R^n$), whereas $f$ itself might not have
such a nice representation.  

{\bf Methodology.} To prove our main result, we proceed in three main
steps. 

1. We first define an infinite dimensional linear programming
problem on an appropriate space of measures, whose optimal value is
the global minimum of $f$ on the set $V$. 

2. We then prove a crucial result, namely that there is {\it no} duality gap between this
linear programming problem and its dual. The approach is similar but
different from that taken in \cite{lasserre} when $V\equiv\R^n$.
Indeed, the approach in 
\cite{lasserre} does {\it not} work when $V\not\equiv\R^n$. Here, we use the
important fact that the polynomial $\theta_r$ is a {\it moment function}.
And so, if a set of probability measures $\Pi$ satisfies
$\sup_{\mu\in\Pi} \int \theta_rd\mu<\infty$, it is {\it tight},
and therefore, by Prohorov's theorem, relatively compact. This latter
intermediate result is crucial for our purpose.

3. In the final step, we use our recent result \cite{lasserre} which states
that if a polynomial $h\in\ring$ is nonnegative on $\R^n$, then
$h+\epsilon\theta_r$ ($\epsilon >0$) is a sum of squares, provided that $r$ is
sufficiently large.

The paper in organized as follows. After introducing the notation and
definitions in \S\ref{notation}, some preliminary results
are stated in  \S\ref{preliminaries}, whereas our main result is stated
and discussed in \S\ref{mainsection}.
For clarity of exposition, most proofs are postponed in 
\S\ref{proofs}, and some auxiliary results are stated in an Appendix;
in particular, duality results for linear programming in 
infinite-dimensional spaces are briefly reviewed.

\section{Notation and definitions}
\label{notation}

Let $\R_+\subset\R$ denote the cone of nonnegative real numbers.
For a real symmetric matrix $A$, the notation $A\succeq0$
(resp. $A\succ0$) stands for
$A$ positive semidefinite (resp. positive definite).
The sup-norm $\sup_j\vert x_j\vert$ of a 
vector $x\in\R^n$, is denoted by $\Vert x\Vert_\infty$.
Let $\ring$ be the ring of real polynomials, and let

\begin{equation}
\label{aa4}
v_r(x):=(1,x_1,x_2,\dots x_n,x_1^2,x_1x_2,\dots,x_1x_n,x_2^2,x_2x_3,\dots,x_n^r)
\end{equation}
be the canonical basis for the $\R$-vector space 
$\mathcal{A}_r$ of real polynomials of degree at most $r$,
and let $s(r)$ be its dimension.
Similarly, $v_\infty(x)$ denotes the canonical basis of $\ring$ as a $\R$-vector
space, denoted $\mathcal{A}$. So a vector in $\mathcal{A}$ has always
{\it finitely} many zeros.

Therefore, a polynomial $p\in\mathcal{A}_r$ is written
\[x\mapsto p(x)\,=\,\sum_{\alpha}p_\alpha
x^\alpha\,=\,\langle \rp,v_r(x)\rangle,\hspace{1cm}x\in\R^n,\]
(where $x^{\alpha}=x_1^{\alpha_1} x_2^{\alpha_2} \dots
x_n^{\alpha_n}$) for some vector 
$\rp=\{ p_\alpha\}\in\R^{s(r)}$, the vector of coefficients of
$p$ in the basis (\ref{aa4}). 

Extending $\rp$ with zeros, we can also
consider $\rp$ as a vector indexed in the basis $v_\infty(x)$
(i.e. $\rp\in\mathcal{A}$). 
If we equip $\mathcal{A}$ with the usual scalar product $\langle
.,.\rangle$ of vectors, then 
for every $p\in\mathcal{A}$,
\[p(x)\,=\,\sum_{\alpha\in\N^n}p_\alpha
x^\alpha\,=\,\langle \rp,v_\infty(x)\rangle,\hspace{1cm}x\in\R^n.\]

Given a sequence $\y=\{y_\alpha\}$ indexed in the basis $v_\infty(x)$,
let $L_\y:\mathcal{A}\to\,\R$ be the linear functional
\begin{equation}
\label{ly}
p\mapsto L_\y(p)\,:=\,\sum_{\alpha\in\N^n}\,p_\alpha
y_\alpha\,=\,\langle \rp ,\y\rangle.
\end{equation}

Given a sequence $\y=\{y_\alpha\}$ indexed in the basis $v_\infty(x)$, the {\it moment}
matrix $M_r(\y)\in\R^{s(r)\times s(r)}$ with
rows and columns indexed in the basis $v_r(x)$ in (\ref{aa4}),
satisfies
\begin{equation}
\label{mommatrix}
\left[M_r(\y)(1,j)\,=\,y_\alpha \:\mbox{ and
}\:M_r(y)(i,1)\,=\,y_\beta\right]\,\Rightarrow\,
M_r(y)(i,j)\,=\,y_{\alpha+\beta}.
\end{equation}
For instance, with $n=2$, 
\[M_2(\y)\,=\,\left[\begin{array}{cccccc}
y_{00}&y_{10}&y_{01} &y_{20}&y_{11}&y_{02}\\
y_{10}&y_{20}&y_{11}&y_{30}&y_{21}&y_{12}\\
y_{01}&y_{11}&y_{02}&y_{21}&y_{12}&y_{03}\\
y_{20}&y_{30}&y_{21}&y_{40}&y_{31}&y_{22}\\
y_{11}&y_{21}&y_{12}&y_{31}&y_{22}&y_{13}\\
y_{02}&y_{12}&y_{03}&y_{22}&y_{13}&y_{04}
\end{array}\right].\]
A sequence $\y=\{y_\alpha\}$ has a {\it representing} measure $\mu_\y$
if
\begin{equation}
\label{momentseq}
y_\alpha\,=\,\int_{\R^n} x^\alpha\,d\mu_\y,\qquad \forall\,\alpha\in\N^n.
\end{equation} 
In this case one also says that $\y$ is a {\it moment sequence}. In addition, if
$\mu_\y$ is unique then $\y$ is said to be a {\it determinate} moment sequence.

The matrix $M_r(\y)$ defines a bilinear form $\langle .,.\rangle_\y$ 
on $\mathcal{A}_r$, by
\[\langle q,p\rangle_\y\,:=\,\langle
\rq,M_r(\y)\rp\rangle\,=\,L_\y(qp),\hspace{0.3cm}q,p\in\mathcal{A}_r,\]
and if $\y$ has a {\it representing} measure $\mu_\y$, then
\begin{equation}
\label{moment1}
L_\y(q^2)\,=\,\langle \rq,M_r(\y)\rq\rangle\,=\,\int_{\R ^n} q(x)^2\,\mu_\y(dx)\,\geq\,0,
\quad\forall \,q\in\mathcal{A}_r, 
\end{equation}
so that $M_r(\y)$ is positive semidefinite, i.e., $M_r(\y)\succeq 0$.

\section{Preliminaries}
\label{preliminaries}

Let $V\subset\R^n$ be the real algebraic set defined in (\ref{setv}), and 
let $B_M$ be the closed ball
\begin{equation}
\label{ball}
B_M\,=\,\{x\in\R^n\,\vert\quad \Vert x\Vert_\infty \,\leq \,M\}.
\end{equation}

\begin{prop}
\label{prop0}
Let $f\in\ring$ be such that $-\infty<f^*:=\inf_{x\in V}f(x)$. 
Then, for every $\epsilon >0$, there is
some $M_\epsilon\in\N$ such that 
\[f^*_M\,:=\,\inf\:\{f(x)\:\vert\: x\in B_M\cap V\}
\,<\,f^*+\epsilon,\qquad \forall M\,\geq \,M_\epsilon.\]
Equivalently, $f^*_M\downarrow f^*$ as $M\to\infty$.
\end{prop}

\begin{proof}
Suppose it is false. That is, there is some $\epsilon_0>0$ and an
infinite sequence sequence $\{M_k\}\subset\N$, with $M_k\to\infty$, such that 
$f^*_{M_k}\geq f^*+\epsilon_0$ for all $k$. But let $x_0\in V$ be such
that $f(x_0)<f^*+\epsilon_0$. With any $M_k\geq \Vert x_0\Vert_\infty$,
one obtains the contradiction $f^*+\epsilon_0\leq f^*_{M_k}\leq f(x_0)<f^*+\epsilon_0$.
\end{proof}

For every $r\in\N$, let $\theta_r\in\ring$ be the polynomial
\begin{equation}
\label{theta}
x\,\mapsto\,\theta_r(x)\,:=\,\sum_{k=0}^r\sum_{i=1}^n\frac{x^{2k}_i}{k{\rm
!}},\qquad x\in\R ^n,
\end{equation}
and notice that $n\leq\theta_r(x)\leq \sum_{i=1}^n\e ^{x^2_i}=:\theta_\infty
(x)$, for all $x\in\R^n$. Moreover, $\theta_r$ is a {\it moment
function}, as it satisfies 
\begin{equation}
\label{momfunc}
\lim_{M\to\infty}\,\inf_{x\in B_M^c} \theta_r(x)\,=\,+\infty,
\end{equation}
where $B_M^c$ denotes the complement of $B_M$ in $\R^n$;
see e.g. Hernandez-Lerma and Lasserre \cite[p. 10]{hernand}.

Next, with $V$ as in (\ref{setv}), introduce the following optimization problems.
\begin{equation}
\label{a0}
\P:\qquad f^*\,:=\,\inf_{x\in V}\,f(x),
\end{equation}
and for $0<M\in\N$, $r\in\N\cup\{\infty\}$,
\begin{equation}
\label{a1}
\mathcal{P}_M^r:\left\{ \begin{array}{lll}
&\displaystyle{\inf_{\mu}}
\int f\,d\mu &\\
\mbox{s.t.}&\int g^2_j\,d\mu&\leq0,\quad j=1,\ldots,m\\
&\int\,\theta_r\,d\mu&\leq n\e ^{M^2}\\
&\mu\in\mathcal{P}(\R^n),&\end{array}\right.
\end{equation}
where $\mathcal{P}(\R^n)$ is the space of probability measures on $\R^n$
(with $\b$ its associated Borel $\sigma$-algebra).
The respective optimal values of $\P$ and $\mathcal{P}_M^r$ are denoted $\inf\P=f^*$ and
$\inf\mathcal{P}_M^r$, or $\min\P$ and 
$\min\mathcal{P}_M^r$ if the minimum is attained (in which case, the
problem is said to be solvable).

\begin{prop}
\label{prop1}
Let $f\in\ring$, and let $\P$ and $\mP^r_M$ be as in
(\ref{a0}) and (\ref{a1}) respectively. Assume that $f^*>-\infty$.
Then,  for every $r\in\N\cup \{\infty\}$, 
$\inf\mathcal{P}_M^r\downarrow f^*$ as $M\to\infty$. If $f$ has a
global minimizer $x^*\in V$, then $\min\mathcal{P}_M^r=f^*$ whenever
$M\geq \Vert x^*\Vert_\infty$.
\end{prop}

\begin{proof}
When $M$ is sufficiently large, $B_M\cap V\neq\emptyset$, and
so, $\mP^r_M$ is consistent, and $\inf\mP^r_M<\infty$.
Let $\mu\in\mP(\R^n)$ be admissible for $\mP_M^r$.
From $\int g_j^2\,d\mu\leq 0$ for all $j=1,\ldots,m$, it follows that 
$g_j(x)^2=0$ for $\mu$-almost all $x\in\R^n$, $j=1,\ldots,m$, 
That is, for every $j=1,\ldots,m$, there exists a set $A_j\in\mathcal{B}$ such that 
$\mu(A^c_j)=0$ and $g_j(x)=0$ for all $x\in A_j$. Take
$A=\cap_jA_j\in\mathcal{B}$ so that $\mu( A ^c)=0$, and 
for all $x\in A$, $g_j(x)=0$ for all $j=1,\ldots,m$. 
Therefore, $A\subset V$, and
as $\mu(A ^c)=0$,
\[\int_{\R^n}f\,d\mu\,=\,\int_{A}f\,d\mu\,\geq \,f^*\quad\mbox{because $f\geq f^*$ on
$A \subset V$,}\]
which proves $\inf\mathcal{P}_M^r\geq f^*$.

As $V$ is closed and $B_M$ is closed and bounded, 
the set $B_M\cap V$ is compact and so, with $f^*_M$ as
in Proposition \ref{prop0}, there is some
$\hat{x}\in B_M\cap V$ such that $f(\hat{x})=f^*_M$. In addition
let $\mu\in\mP(\R^n)$ be the Dirac probability measure at the
point $\hat{x}$. As $\Vert\hat{x}\Vert_\infty\leq M$,
\[\int\,\theta_r\,d\mu\,=\,
\theta_r(\hat{x})\,\leq\, n\e ^{M^2}.\]
Moreover, as $\hat{x}\in V$, $g_j(\hat{x})=0$, for all $j=1,\ldots,m$,
and so
\[\int\,g^2_j\,d\mu\,=\,
g_j(\hat{x})^2\,=\,0,\quad j=1,\ldots,m,\]
so that $\mu$ is an admissible solution of $\mathcal{P}_M^r$
with value $\int f\,d\mu=f(\hat{x})=f^*_M$, which proves that $\inf\mathcal{P}_M^r
\leq f^*_M$. This latter fact, combined with 
Proposition \ref{prop0} and with $f^*\leq\inf\mathcal{P}_M^r$, implies
$\inf\mathcal{P}_M^r\downarrow f^*$ as $M\to\infty$, the desired
result. The final statement is immediate by taking as feasible
solution for $\mathcal{P}_M^r$, the Dirac probability measure at the
point $x^*\in B_M\cap V$
(with $M\geq \Vert x^*\Vert_\infty$). As its value is now $f^*$, it is
also optimal, and so, $\mathcal{P}_M^r$ is solvable with optimal value $\min\mathcal{P}_M^r=f^*$.
\end{proof}

Consider now, the following optimization problem $\mathcal{Q}^r_M$, the dual problem of
$\mathcal{P}^r_M$, i.e.,
\begin{equation}
\label{duala1}
\mathcal{Q}^r_M:\quad\begin{array}{lll}
&\displaystyle{\max_{\lambda,\delta,\gamma}} \quad\gamma -n\delta\e ^{M^2}&\\
\mbox{s.t.}& f+\delta\theta_r+\sum_{j=1}^m\lambda_jg_j^2&\geq \gamma\\
&\gamma\in\R,\delta\in\R_+,\lambda\in\R^m_+,&
\end{array} 
\end{equation}
with optimal value denoted by $\sup\mathcal{Q}^r_M$.
Indeed, $\mathcal{Q}^r_M$ is a dual of $\mP^r_M$ because {\it weak duality}
holds. To see this, consider any two feasible solutions $\mu\in\mathcal{P}(\R^n)$ and 
$(\lambda,\delta,\gamma)\in\R^m_+\times\R_+\times\R$, 
of $\mP^r_M$ and $\mathcal{Q}^r_M$,
respectively.
Then, integrating both sides of the inequality in  $\mathcal{Q}^r_M$
with respect to $\mu$, yields
\[\int fd\mu +\delta\int \theta_r\,d\mu+
\sum_{j=1}^m\lambda_j\int \,g_j^2\,d\mu\geq \gamma,\]
and so, using that $\mu$ is feasible for $\mP^r_M$,
\[\int fd\mu \,\geq\,
\gamma -\delta n\e ^{M^2}.\]
Hence, the value of
any feasible solution of $\mathcal{Q}^r_M$ is always smaller than the
value of any feasible solution of $\mP^r_M$, i.e., weak duality holds.

In fact we can get the more important and crucial following result.

\begin{theorem}
\label{dualmom}

Let $M$ be large enough so that $B_M\cap V\neq\emptyset$.
Let $f\in\ring$, and let $r_0>\max[{\rm deg}\, f,{\rm deg}\,g_j]$.
Then, for every $r\geq r_0$, 
$\mathcal{P}^r_M$ is solvable, and
there is no duality gap between
$\mathcal{P}^r_M$ and its dual $\mathcal{Q}^r_M$. That is,
$\sup\mathcal{Q}^r_M=\min\mathcal{P}^r_M$.
\end{theorem}
For a proof see \S\ref{proofdualmom}.
We finally end up this section by re-stating a result proved in
\cite{lasserre}, which, together 
with Theorem \ref{dualmom}, will be crucial to prove our main result.

\begin{theorem}[\cite{lasserre}]
\label{uncons}
Let $f\in\ring$ be nonnegative.
Then for every $\epsilon >0$, there is some $r(\epsilon)\in\N$
such that, 
\begin{equation}
\label{uncons-2}
f_{\epsilon r(\epsilon)}\:(=f+\epsilon\theta_{r(\epsilon)})
\qquad\mbox{is a sum of squares,}
\end{equation}
and so is $f_{\epsilon r}$, for all $r\geq r(\epsilon)$.
\end{theorem}

\section{Main result}
\label{mainsection}

Recall that for given $(\epsilon,r)\in\R\times\N$,
$f_{\epsilon r}=f+\epsilon\theta_r$,
with $\theta_r\in\ring$ being the polynomial defined in (\ref{theta}).
We now state our main result:

\begin{theorem}
\label{thmain}
Let $V\subset\R^n$ be as in (\ref{setv}), and let $f\in\ring$ be
nonnegative on $V$. Then, for every $\epsilon>0$, there exists
$r(\epsilon)\in\N$ and nonnegative scalars $\{\lambda_j\}_{j=1}^m$,
such that, for all $r\geq r(\epsilon)$,
\begin{equation}
\label{main-1}
f_{\epsilon r}\,=\,q-\sum_{j=1}^m\lambda_j\,g_j^2,
\end{equation}
for some s.o.s. polynomial $q\in\ring$.
In addition, $\Vert f-f_{\epsilon r}\Vert_1\to 0$, as
$\epsilon\downarrow 0$.
\end{theorem}
For a proof see \S\ref{proofmain}.

\begin{remark}
\label{remark-1}
(i) Observe that (\ref{main-1}) is an obvious certificate of
positivity of $f_{\epsilon r}$ on the algebraic set $V$, because everywhere
on $V$, $f_{\epsilon r}$ coincides with the s.o.s. polynomial $q$. Therefore, when $f$
is nonnegative on $V$, one obtains with {\it no} assumption on the
algebraic set $V$, a certificate of positivity for any approximation
$f_{\epsilon r}$ of $f$ (with $r\geq r(\epsilon)$), whereas $f$ itself might not have such a
representation. In other words, the $(\epsilon,r)$--perturbation $f_{\epsilon r}$ of $f$, has
a {\it regularization} effect on $f$ as it permits to derive nice
representations.

(ii) From the proof of Theorem \ref{thmain}, instead of the
representation (\ref{main-1}), one may also provide
the alternative representation 
\[f_{\epsilon r}\,=\,q-\lambda\,\sum_{j=1}^m g_j^2,\]
for some s.o.s. polynomial $q$, and some 
(single) nonnegative scalar $\lambda$ (instead of $m$
nonnegative scalars in (\ref{main-1})).
\end{remark}

\subsection{The case of a semi-algebraic set}

We now consider the representation of polynomials, nonnegative on a
semi algebraic set $\K\subset\R^n$, defined as,
\begin{equation}
\label{setk}
\K\,:=\,\{x\in\R^n\,\vert\quad g_j(x)\,\geq\,0,\quad j=1,\ldots,m\},
\end{equation}
for some family $\{g_j\}_{j=1}^m\subset \ring$.

One may apply the machinery developed previously for algebraic sets, because
the semi-algebraic set $\K$ may be viewed as the projection on $\R^n$,
of an algebraic set in
$\R^{n+m}$. Indeed, let $V\subset\R^{n+m}$ be the algebraic set defined as
\begin{equation}
\label{newsetv}
V\,:=\,\{(x,z)\in\R^n\times\R^m\,\vert\quad g_j(x)-z_j^2\,=\,0,\quad j=1,\ldots,m\}.
\end{equation}
Then every $x\in\K$ is associated with the point
$(x,\sqrt{g_1(x)},\ldots,\sqrt{g_m(x)})\in V$.

Let $\R[z]:=\R[z_1,\ldots,z_m]$, and $\R[x,z]:=\R[x_1,\ldots x_n,z_1,\ldots,z_m]$, and 
for every $r\in\N$, let
$\varphi_r\in\R[z]$ be the polynomial
\begin{equation}
\label{varphi}
z\,\mapsto\,\varphi_r(z)\;=\,\sum_{k=0}^r\sum_{j=1}^m\frac{z_j^{2k}}{k{\rm
!}}.
\end{equation}

We then get :
\begin{corollary}
\label{corosemialg}
Let $\K$ be as in (\ref{setk}), and
$\theta_r,\varphi_r$ be as in (\ref{theta}) and (\ref{varphi}). 
Let $f\in\ring$ be nonnegative on $\K$. 
Then, for every $\epsilon >0$, there exist 
nonnegative scalars $\{\lambda_j\}_{j=1}^m$ such that,
for all $r$ sufficiently large,
\begin{equation}
\label{semi-1}
f+\epsilon\theta_r+\epsilon\varphi_r\,=\,q_\epsilon-\sum_{j=1}^m
\lambda_j(g_j-z_j^2)^2,
\end{equation}
for some s.o.s. polynomial $q_\epsilon\in\R[x,z]$.

Equivalently,  everywhere on $\K$, the polynomial
\begin{equation}
\label{coro-2}
x\mapsto f(x)+\epsilon\sum_{k=0}^r\sum_{i=1}^n\frac{x^{2k}_i}{k{\rm !}}\,+\,
\epsilon\sum_{k=0}^r\sum_{j=1}^m\frac{g_j(x)^{k}}{k{\rm !}}
\end{equation}
coincides with the nonnegative polynomial
$x\mapsto q_\epsilon(x,\sqrt{g_1(x)},\ldots,\sqrt{g_m(x)})$.
\end{corollary}

So, as for the case of an algebraic set $V\subset\R^n$,
 (\ref{semi-1}) is an obvious certificate of positivity on the
semi-algebraic set $\K$, for the
polynomial $f_{\epsilon r}\in\R[x,z]$
\[f_{\epsilon r}\,:=\,f+\epsilon\theta_r+\epsilon\varphi_r,\]
and in addition, viewing $f$ as an element of $\R[x,z]$, one has
$\Vert f-f_{\epsilon r}\Vert_1\to 0$ as
$\epsilon\downarrow  0$. Notice that {\it no} assumption on $\K$ or on
the $g_j$'s that define $\K$, is needed.

Now, assume that $\K$ is compact and the $g_j$'s that define $\K$,
satisfy Putinar's condition, i.e., (i) there exits some $u\in\ring$ such
that $u$ can be written $u_0+\sum_ju_jg_j$ for some s.o.s. polynomials
$\{u_j\}_{j=0}^m$, and (ii), the level set $\{x\vert u(x)\geq 0\}$ is
compact.

If $f$ is nonnegative on $\K$, then $f+\epsilon\theta_r$ is strictly
positive on $\K$, and therefore, by Putinar's theorem \cite{putinar}
\begin{equation}
\label{put}
f+\epsilon\theta_r\,=\,q_0+\sum_{j=1}^mq_j g_j,
\end{equation}
for some s.o.s. family $\{q_j\}_{j=0}^m$. One may thus either have
Putinar's representation (\ref{put}) in $\R^n$, or 
(\ref{semi-1}) via a {\it lifting} in $\R^{n+m}$.

One may relate (\ref{semi-1}) and (\ref{put}) by
\[q_\epsilon(x,z)\,=\,q_{\epsilon}^1(x)+q_{\epsilon }^2(x,z^2),\]
with
\[x\,\mapsto\,q_{\epsilon}^1(x)\,:=\,q_0(x)+\sum_{j=1}^m
\left(q_j(x)g_j(x)+\lambda_jg_j(x)^2\right),\]
and
\[(x,z)\,\mapsto\,q_{\epsilon}^2(x,z^2)\,:=\,
\epsilon\varphi_r(z)+\sum_{j=1}^m\lambda_jz_j^4-2g_j(x)z_j^2.\]

\subsection{Computational implications}

The results of the previous section can be applied to compute (or at
least approximate) the global minimum of $f$ on $V$.
Indeed, with $\epsilon>0$ fixed, and $2r\geq\max[{\rm deg}\,f, {\rm deg}\,g_j^2]$,
consider the convex optimization problem
\begin{equation}
\label{primalagain}
\Q_{\epsilon r}\left\{\begin{array}{llcl}
&\displaystyle{\min_\y\:L_\y(f_{\epsilon r})},&&\\
{\rm s.t.}&M_{r}(\y)&\succeq &0\\
&L_\y(g_j^2)&\leq &0,\quad j=1,\ldots,m\\
&y_{0}&=&1,
\end{array}\right.
\end{equation}
where $\theta_r$ is as in (\ref{theta}), $L_\y$ and $M_r(\y)$ are 
the linear functional and the moment matrix associated with a sequence
$\y$ indexed in the basis (\ref{aa4}); see (\ref{ly}) and
(\ref{mommatrix}) in \S\ref{notation}.

$\Q_{\epsilon r}$ is called a semidefinite programming (SDP) problem,
and its associated {\it dual} SDP problem
reads
\begin{equation}
\label{dualagain}
\Q^*_{\epsilon r}\left\{\begin{array}{lll}
&\displaystyle{\max_{\lambda,\gamma,q} \quad\gamma } &\\
{\rm s.t.}&f_{\epsilon r}-\gamma&=
\displaystyle{q-\sum_{j=1}^m\lambda_j g_j^2}\\
&\lambda\in\R^m, &\lambda\geq0,\\
&q\in\ring, &\mbox{$q$ s.o.s. of degree }\leq 2r.
\end{array}\right.
\end{equation}
Their optimal values are denoted $\inf\Q_{\epsilon r}$ and
$\sup\Q^*_{\epsilon r}$, respectively
(or $\min\Q_{\epsilon r}$ and $\max\Q^*_{\epsilon r}$ 
if the optimum is attained, in which case the problems
are said to be solvable). 
Both problems $\Q_{\epsilon r}$ and its dual $\Q^*_{\epsilon r}$ are
nice convex optimization problems that, in principle, can be solved
efficiently by standard software packages.
For more details on SDP theory, the interested
reader is referred to the survey paper \cite{boyd}.

That weak duality holds between $\Q_{\epsilon r}$ and $\Q^*_{\epsilon
r}$ is straightforward. Let $\y=\{y_\alpha\}$ and
$(\lambda,\gamma,q)\in\R^m_+\times\R\times \ring$ be feasible
solutions of $\Q_{\epsilon r}$ and $\Q^*_{\epsilon r}$,
respectively. Then, by linearity of $L_\y$,
\begin{eqnarray*}
L_\y(f_{\epsilon r}) -\gamma&=&L_\y(f_{\epsilon r} -\gamma)\\
&=&L_\y(q-\sum_{j=1}^m\lambda_j
g_j^2)\,=\,L_\y(q)-\sum_{j=1}^m\lambda_j L_\y(g_j^2)\\
&\geq& L_\y(q) \quad\mbox{[because $L_\y(g_j^2)\leq 0$ for all $j=1,\ldots,m$]}\\
&\geq& 0\quad\mbox{[because $q$ is s.o.s. and $M_r(\y)\succeq 0\,$; see (\ref{moment1}).]}
\end{eqnarray*}
Therefore, $L_\y(f_{\epsilon r})\geq \gamma$, the desired conclusion.
Moreover, $\Q_{\epsilon r}$ is an obvious relaxation of the perturbed
problem 
\[\P_{\epsilon r}:\quad f^*_{\epsilon r}:=\min_{x}\:\{f_{\epsilon r}\,\vert\:x\in V\}.\]
Indeed, let $x\in V$ and let $\y:=v_{2r}(x)$ (see (\ref{aa4})), i.e., $\y$
is the vector of moments (up to order $2r$) of the Dirac measure at
$x\in V$. Then, $\y$ is feasible for $\Q_{\epsilon r}$ because
$y_0=1$, $M_r(\y)\succeq 0$, and $L_\y(g_j^2)=g_j(x)^2=0$, for all $j=1,\ldots ,m$. Similarly,
$L_\y(f_{\epsilon r})=f_{\epsilon r}(x)$. Therefore, $\inf\Q_{\epsilon
r}\leq f^*_{\epsilon r}$.

\begin{theorem}
\label{thcomp}
Let $V\subset\R^n$ be as in (\ref{setv}), and $\theta_r$ as in
(\ref{theta}). Assume that $f$ has a global minimizer $x^*\in V$ with
$f(x^*)=f^*$. Let $\epsilon>0$ be fixed.
Then
\begin{equation}
\label{thcomp-1}
f^*\,\leq\,\sup\Q^*_{\epsilon r}\,\leq\,\inf\Q_{\epsilon r}\,\leq\,f^*+\epsilon\theta_r(x^*)\,\leq\,
f^*+\epsilon\sum_{i=1}^n\e ^{(x^*_i)^2},
\end{equation}
provided that $r$ is sufficiently large.
\end{theorem}
\begin{proof}
Observe that the polynomial $f-f^*$ is nonnegative on $V$. Therefore,
by Theorem \ref{thmain}, for every $\epsilon$ there exists
$r(\epsilon)\in\N$ and $\lambda(\epsilon)\in\R^m_+$, such that 
\[f-f^*+\epsilon\theta_r+\sum_{j=1}^m\lambda_j(\epsilon)g_j^2\,=\,q_\epsilon,\]
for some s.o.s. polynomial $q_\epsilon\in\ring$. 
But this shows that
$(\lambda(\epsilon),f^*,q_\epsilon)\in\R^m_+\times\R\times\ring$ is a feasible
solution of $\Q^*_{\epsilon r}$ as soon as $r\geq r(\epsilon)$, in which case,
$\sup\Q^*_{\epsilon r}\geq f^*$. Moreover, we have
seen that $\inf\Q_{\epsilon r}\leq f_{\epsilon r}(x)$ for
any feasible solution $x\in V$. In particular,
$\inf\Q_{\epsilon r}\leq f^*+\epsilon\theta_r(x^*)$, from which (\ref{thcomp-1}) follows.
\end{proof}

Theorem \ref{thcomp} has a nice feature. Suppose that one
knows some bound $\rho$ on the norm $\Vert x^*\Vert_\infty$ of a global
minimizer of $f$ on $V$. Then, one may fix \`a priori the error bound
$\eta$ on $\vert\inf\Q_{\epsilon r}-f^*\vert$. Indeed, let $\eta$ be fixed, and 
fix $\epsilon>0$ such that
$\epsilon\leq \eta(n\e ^{\rho ^2})^{-1}$. By Theorem \ref{thcomp}, one
has $f^*\leq \inf\Q_{\epsilon r}\leq f^*+\eta$, provided that $r$ is large enough.

The same approach works to approximate the global minimum of a 
polynomial $f$ on a semi-algebraic set $\K$, as defined in (\ref{setk}). 
In view of Corollary \ref{corosemialg}, and via a lifting in $\R^{n+m}$,
one is reduced to the case of a real algebraic set $V\subset\R^{n+m}$, so that
Theorem \ref{thcomp} still applies. It is important to emphasize that
one requires {\it no} assumption on $\K$, or on the $g_j$'s that define
$\K$. This is to be compared with previous SDP-relaxation techniques
developed in e.g. \cite{lasserre1,lasserre2,lasserre3,parrilo2,markus},
where the set $\K$ is supposed to be
compact, and with an additional assumption on the $g_j$'s to ensure that
Putinar's representation \cite{putinar} holds.

\section{Proofs}
\label{proofs}

\subsection{Proof of Theorem \ref{dualmom}}
\label{proofdualmom}
To prove the absence of a duality gap, we first rewrite
$\mP^r_M$ (resp. $\mathcal{Q}^r_M$) as a linear program in
(standard) form
\[\min_x \{\ip{x}{c}\,\vert \quad Gx=b, x\in C\},\quad \mbox{(resp. }
\max_w \{\ip{w}{b}\,\vert \quad c-G^*w\in C ^*\}\mbox{)},\]
on appropriate dual pairs of vector spaces, with associated
convex cone $C$ (and its dual $C^*$), and associated linear map
$G$ (and its adjoint $G^*$). Then, we will prove that $G$ is
continuous, and the set 
$D:=\{(Gx,\ip{x}{c})\,\vert\,x\in C\}$
is closed, in some appropriate weak topology. This permits us to conclude by invoking
standard results in infinite-dimensional linear programming, 
that one may find in e.g. Anderson and Nash
\cite{nash}. For a brief account see \S\ref{LP}, and for more details,
see e.g. Robertson and Robertson \cite{robertson},
and Anderson and Nash \cite{nash}.

Let $\theta_r$ be as in (\ref{theta}), and let
$M(\R^n)$ be the $\R$-vector space of finite signed Borel measures $\mu$ on $\R^n$,
such that $\int \theta_r\,d\vert\mu\vert<\infty$ (where
$\vert\mu\vert$ denotes the total variation of $\mu$). 
Similarly, let $H^r$ be the $\R$-vector space of continuous functions $h:\R^n\to\R$, such that
$\sup_{x\in\R^n}\vert h(x)\vert/\theta_r(x)<\infty$.
With the bilinear form $\ip{.}{.}:\:M(\R^n)\times H^r$, defined as
\[(\mu,h)\,\mapsto\,\ip{\mu}{h}\,=\,\int h\,d\mu,\qquad (\mu,h)\in
M(\R^n)\times H^r,\]
$(M(\R^n),H^r)$ forms a {\it dual} pair of vector spaces 
(See \S\ref{LP}.) Introduce the dual pair of vector spaces $(\X,\Y)$,
\[\X\,:=\,M(\R^n)\times \R^m\times\R,\quad
\Y\,:=\,H^r\times\R^m\times\R,\]
and $(\Z,\W)$
\[\Z\,:=\,\R^m\times\R\times\R,\quad
\W\,:=\,\R^m\times\R\times \R.\]
Recall that $2r>{\rm deg}\,g_j^2$, for all $j=1,\ldots,m$,
and let $G:\X\to\Z$ be the linear map 
\[(\mu,u,v)\,\mapsto\,G(\mu,u,v)\,:=\,\left[\begin{array}{l}\ip{\mu}{g_1^2}+u_1\\
\ldots\\ \ip{\mu}{g_m^2}+u_m\\
\ip{\mu}{\theta_r}+v\\
\ip{\mu}{1}\end{array}\right],\]
with associated {\it adjoint} linear map $G^*:\W\to\Y$
\[(\lambda,\delta,\gamma)\,\mapsto\,G^*(\lambda,\delta,\gamma)\,
:=\,\left[\begin{array}{l}\sum_{j=1}^m\lambda_jg_j^2+\delta\theta_r+\gamma\\
\lambda\\\delta
\end{array}\right],\]
Next, let $M(\R^n)_+\subset M(\R^n)$ be the convex cone of 
nonnegative finite Borel measures on $\R^n$, so that the set
$C:= M(\R^n)_+\times \R^m_+\times \R_+\subset \X$
is a convex cone in $\X$. If $H^r_+$ denotes the nonnegative functions
of $H^r$, then 
\[C^*\,=\,H^r_+\times\R^m_+\times\R_+\subset\,\Y.\]
is the {\it dual} cone of $C$ in $\Y$.

As $2r>\max[{\rm deg }f,{\rm deg}\,g_j^2]$ it follows that
$f\in H^r$ and  $g_j^2\in H^r$, for all $j=1,\ldots,m$.
Then, by introducing {\it slack} variables $u\in\R^m_+,v\in\R_+$, rewrite 
the infinite-dimensional linear program $\mP^r_M$ defined in (\ref{a1}), in
equality form, that is,
\begin{equation}
\label{aa1}
\mathcal{P}_M^r:\left\{ \begin{array}{ll}
&\displaystyle{\inf_{\mu,u,v}}\:
\ip{(\mu,u,v)}{(f,0,0)} \\
\mbox{s.t.}&G(\mu,u,v)\,=\,\left[\begin{array}{l}0\\n\e ^{M^2}\\1\end{array}\right]\\
&(\mu,u,v)\in C.\end{array}\right.
\end{equation}
The LP dual $(\mP^r_M)^*$ of $\mP^r_M$ now reads
\begin{equation}
\label{aa2}
(\mP^r_M)^*:\left\{ \begin{array}{ll}
&\displaystyle{\sup_{\lambda,\delta,\gamma}}\:
\ip{(\lambda,\delta,\gamma)}{(0,n\e ^{M^2},1)}\\
\mbox{s.t.}&(f,0,0)-G ^*(\lambda,\delta,\gamma)\in C^*.
\end{array}\right.
\end{equation}
Hence, every feasible solution $(\lambda,\delta,\gamma)$ of
$(\mP^r_M)^*$ satisfies
\begin{equation}
\label{feasible}
f-\sum_{j=1}^m\lambda_j\,g_j^2-\delta\,\theta_r-\gamma\geq0;\quad
\lambda,\delta\leq0.
\end{equation}
As $\lambda,\delta\leq 0$ in (\ref{aa2}), one may see that the two formulations
(\ref{aa2}) and (\ref{duala1}) are identical, i.e., $\mathcal{Q}^r_M=(\mP^r_M)^*$.

As $2r>\max[{\rm deg }f,{\rm deg}\,g_j^2]$, it follows that
$f-\sum_{j=1}^m\lambda_j\,g_j^2-\delta\,\theta_r-\gamma\in H^r$,
for all $(\lambda,\delta,\gamma)\in\W$. Therefore,
$G^*(\W)\subset\Y$, and so, by Proposition \ref{proplp},
the linear map $G$ is weakly continuous (i.e. is continuous with respect
to the weak topologies $\sigma(\X,\Y)$ and $\sigma(\Z,\W)$). 

We next prove that the set $D\subset\Z\times\R$, defined as
\begin{equation}
\label{setD}
D\,:=\,\{(G(\mu,u,v),\ip{(\mu,u,v)}{(f,0,0)})\,\vert\quad
(\mu,u,v)\in C\},
\end{equation}
is {\it weakly closed}.

For some directed set $(A,\geq)$, 
let  $\{(\mu_\beta,u_\beta,v_\beta)\}_{\beta\in A}$ be a net in $C$,
such that 
\[(G(\mu_\beta,u_\beta,v_\beta),
\ip{(\mu_\beta,u_\beta,v_\beta)}{(f,0,0)})\,\to\,((a,b,c),d),\]
weakly, for some element $((a,b,c) ,d)\in\Z\times\R$. In particular
\[\mu_\beta(\R^n)\,\to\,c;\quad \ip{\mu_\beta}{\theta_r}+v_\beta\,\to\,b;\quad 
\ip{\mu_\beta}{g_j^2}+(u_\beta)_j\,\to\,a_j,\:j=1,\ldots,m,\]
and $\ip{\mu_\beta}{f}\to d$.
As $(\mu_\beta,u_\beta,v_\beta)\in C$, 
and $\theta_r,g_j^2\geq 0$, it follows
immediately that $a,b,c\geq 0$.
We need to consider the two cases $c=0$ and $c>0$.

Case $c=0$. From $\mu_\beta(\R^n)\to\,c$, it follows that $\mu_\beta\to\mu:=0$
in the total variation norm. 
But in this case, observe that  $G(\mu,a,b)=(a,b,c)$. It remains to prove
that we also have $\ip{\mu_\beta}{f}\to d=0$, in which case,
$(G(\mu,a,b),\ip{\mu}{f})=((a,b,c),d)$, as desired.

Recall that $r\geq \mbox{deg}f$. 
Denote by $\{y_\alpha(\beta)\}_{\vert\alpha\vert\leq 2r}$ the sequence
of moments of the measure $\mu_\beta$, i.e.,
\[y_\alpha(\beta)\,=\,\int x^\alpha\,d\mu_\beta,\quad\alpha\in\N^n,\quad
\vert\alpha\vert\leq 2r.\]
In particular, $y_0(\beta)=\mu_\beta(\R^n)$.
From $\ip{\mu_\beta}{\theta_r}+v_\beta\to\,b$,
there is some $\beta_0\in A$, such that 
$\ip{\mu_\beta}{\theta_r}\leq 2b$ for all $\beta\geq \beta_0$.
But this implies that 
\[y_{2k}(i,\beta)\,:=\,\int x_i^{2k}\,d\mu_\beta\leq 2r{\rm !}b,\quad k\leq r,\quad
i=1,\ldots,n.\]
By Lemma \ref{newlemma}, it follows that $y_{2\alpha}(\beta)\leq 2br{\rm !}$ for
all $\alpha\in\N^n$ with $\vert\alpha\vert\leq r$, and
$\vert y_\alpha(\beta)\vert \leq \sqrt{2y_0(\beta)\,br{\rm !}}$ for all
$\vert\alpha\vert\leq r$. But then, as $y_0(\beta)=\mu_\beta(\R^n)\to
c=0$, we thus obtain $y_\alpha(\beta)\to 0$ for all
$\vert\alpha\vert\leq r$. Therefore,
\[\ip{\mu_\beta}{f}\,=\,\int f\,d\mu_\beta\,=\,\sum_{\vert\alpha\vert\leq
r}f_\alpha \int x^\alpha\,d\mu_\beta\,=\,
\sum_{\vert\alpha\vert\leq r}f_\alpha y_\alpha(\beta)\to 0,\]
the desired result.

Case $c>0$. From $\mu_\beta(\R^n)\to\,c$ and $\ip{\mu_\beta}{\theta_r}+v_\beta\to\,b$,
there is some $\beta_0\in A$, such that $\mu_\beta(\R^n)\leq 2c$ and 
$\ip{\mu_\beta}{\theta_r}\leq 2b$ for all $\beta\geq \beta_0$.
But, as $\theta_r$ is a {\it moment} function, this implies that
the family $\Delta:=\{\nu_\beta:=\mu_\beta/\mu_\beta(\R^n)\}_{\beta\geq  \alpha_0}$ is a
{\it tight} family of probability measures, and as $\Delta$ is a set of probability measures
on a metric space, by Prohorov's theorem,  $\Delta$ is relatively
compact (see \cite[Chap. 1]{hernand} and section \S\ref{aux}).
Therefore, there is some probability measure 
$\nu ^*\in M(\R^n)$, and a sequence $\{n_k\}\subset\Delta$,
such that $\nu_{n_k}$ converges to
$\nu ^*$, for the {{\it weak convergence} of probability measures, i.e.,
\[\ip{\nu_{n_k}}{h}\,\to\,\ip{\nu ^*}{h},\qquad \forall h\in\cb\]
(where $\cb$ denotes the space of bounded continuous functions on $\R^n$);
see e.g. Billingsley \cite{billingsley}.
Hence, with $\mu ^*:=c\,\nu ^*$, we also conclude
\begin{equation}
\label{rweak}
\ip{\mu_{n_k}}{h}\,\to\,\ip{\mu ^*}{h},\qquad \forall h\in\cb.
\end{equation}
Next, as $2r>\max[{\rm deg}\,f,{\rm deg}\,g_j^2]$, 
the functions $f/\theta_{r-1}$ and $g_j^2/\theta_{r-1}$, $j=1,\ldots,m$, are
all in $\cb$. Therefore, using Lemma
\ref{lemmamom}, we obtain
\[\ip{\nu_{n_k}}{f}\,\to\,\ip{\nu ^*}{f},\:\mbox{ and }\:
\ip{\nu_{n_k}}{g_j^2}\,\to\,\ip{\nu ^*}{g_j^2},\:j=1,\ldots,m.\]
And, therefore, 
\[\ip{\mu_{n_k}}{f}\,\to\,\ip{\mu ^*}{f}\,=\,d,\:\mbox{ and }\:
\ip{\mu_{n_k}}{g_j^2}\,\to\,\ip{\mu ^*}{g_j^2},\:j=1,\ldots,m.\]
Finally, from the weak convergence (\ref{rweak}),
and  as $\theta_r$ is continuous and nonnegative,
\[\ip{\mu^*}{\theta_r}\,\leq\,\liminf_{k\to\infty}\:\ip{\mu_{n_k}}{\theta_r}\,\leq\,b,\]
see e.g. \cite[Prop. 1.4.18]{hernand}.

So, let $v:=b-\ip{\mu ^*}{\theta_r}\geq 0$, and
$u_j:=a_j-\ip{\mu ^*}{g_j^2}\geq 0$, $j=1,\ldots,m$, and recalling that
$c=\mu ^*(\R^n)$, we conclude that
$G(\mu ^*,u,v)=(a,b,c)$, and
$\ip{(\mu ^*,u,v)}{(f,0,0)}=d$, which proves that the set $D$ in (\ref{setD}) is weakly closed.

Finally, by Proposition 
\ref{prop1}, $\mP^r_M$ is consistent with finite value as soon as $M$ is large
enough to ensure that $B_M\cap V\neq\emptyset$.
Therefore, one may
invoke Theorem \ref{nash}, and conclude that there is no duality gap
between $\mP^r_M$ and its dual $\mathcal{Q}^r_M$, the desired
result.  $\qed$

\subsection{Proof of Theorem \ref{thmain}}
\label{proofmain}

It suffices to prove the result for the case where  $\inf_{x\in V}f(x)=f ^*>0$.
Indeed, suppose that $f^*=0$. Then with $\epsilon>0$ fixed,
arbitrary, $f^*+n\epsilon>0$ and so, 
suppose that (\ref{main-1}) holds for 
$\hat{f}:=f+n\epsilon$. There is some $r(\epsilon)\in\N$ such that, for all 
$r\geq r(\epsilon)$,
\[\hat{f}\,=\,f+n\epsilon +
\epsilon\,\theta_r \,=\,q_{\epsilon r}-\sum_{j=1}^m\lambda_j g_j^2,\]
for some s.o.s. polynomial $q_{\epsilon r}$, and some nonnegative
scalars $\{\lambda_j\}$.
Equivalently,
\[f+2\epsilon\,\theta_r\,=\,q_{\epsilon r}+
\epsilon\sum_{k=1}^{r}\sum_{j=1}^n
\frac{x_j^{2k}}{k{\rm !}}-\sum_{j=1}^m\lambda_j g_j^2\,=\,
\hat{q}_{\epsilon r}-\sum_{j=1}^m\lambda_j g_j^2,\]
where $\hat{q}_{\epsilon r}$ is a s.o.s. polynomial. Equivalently,
$f_{2\epsilon r}=\hat{q}_{\epsilon r}-\sum_{j=1}^m\lambda_j g_j^2$,
so that (\ref{main-1}) also holds for $f$. Therefore, from now on, we
will assume that $f^*>0$.

So let $\epsilon>0$ (fixed) be such that $f^*-\epsilon>0$, and let
$r\geq r_0$ with $r_0$ as in Theorem \ref{dualmom}. Next, by
Proposition \ref{prop1}, let $M$ be such that 
$f^*\leq \inf\mP^r_M\leq f^*+\epsilon$.
By Theorem \ref{dualmom}, we then have $\sup\mathcal{Q}^r_M\geq
f^*$. So, by considering a maximizing sequence of
$\mathcal{Q}^r_M$, there is some
$(\lambda,\delta,\gamma)\in\R^m_+\times\R_+\times\R$, 
such that
\begin{equation}
\label{pmain0}
0<f^*-\epsilon<\gamma-n\delta\e ^{M^2}\leq f^*+\epsilon;\quad
f+\delta\theta_r+\sum_{j=1}^m\lambda_j\,g_j^2\geq \gamma,
\end{equation}
and so,
\begin{equation}
\label{pmain1}
f-(\gamma
-n\delta\e ^{M^2})+\sum_{j=1}^m\lambda_j\,g_j^2\geq \delta(n\e
^{M^2}-\theta_r).
\end{equation}
By Proposition \ref{prop0}, we may choose $M$ such that
there is some $x_M\in B_{M/2}\cap V$ such that $f(x_M)\leq
f^*+\epsilon$. Evaluating (\ref{pmain1}) at  $x=x_M$ yields
\begin{equation}
\label{pmain2}
2\epsilon\,\geq\,f(x_M)-(\gamma
-n\delta\e ^{M^2})\geq \delta(n\e
^{M^2}-\theta_r(x_M)),
\end{equation}
and so, using $\Vert x_M\Vert_\infty\leq M/2$,
\begin{equation}
\label{pmain3}
2\epsilon\geq \delta n(\e^{M^2}-\e^{M^2/4}),
\end{equation}
which yields $\delta\leq 2\epsilon/n(\e^{M^2}-\e^{M^2/4})$. Therefore,
given $\epsilon>0$, one may pick
$(\lambda,\delta,\gamma)$ in a maximizing sequence of $\mathcal{Q}^r_M$,
in such a way that $\delta\leq \epsilon$.

For such a choice of $(\lambda,\delta,\gamma)$, and in view of (\ref{pmain0}), we have
\[f+\delta\theta_r+\sum_{j=1}^m\lambda_j\,g_j^2\geq
(\gamma-n\delta\e ^{M^2})+n\delta\e ^{M^2}\,\geq\,f^*-\epsilon+n\delta\e ^{M^2}
\geq 0,\]
so that the polynomial $h:=f+\delta\theta_r+\sum_{j=1}^m\lambda_j\,g_j^2$ is nonnegative.

Therefore, invoking Theorem \ref{uncons} proved in Lasserre \cite{lasserre}, there is some
$r(\epsilon)\in\N$ such that, for all $s\geq r(\epsilon)$, the polynomial
$q_\epsilon:=h+\epsilon\theta_s$ is a s.o.s. But then, take $s>\max[r,r(\epsilon)]$ and
observe that
\[\delta\theta_r+\epsilon\theta_s\,=\,(\delta+\epsilon)\theta_s-
\delta\sum_{k=r+1}^s\sum_{j=1}^n\frac{x^2_i}{k{\rm !}},\]
and so
\[q_\epsilon=h+\epsilon\theta_s=f+\sum_{j=1}^m\lambda_jg_j^2+
(\delta+\epsilon)\theta_s-
\delta\sum_{k=r+1}^s\sum_{i=1}^n\frac{x^2_i}{k{\rm
!}},\]
or, equivalently,
\[f+\sum_{j=1}^m\lambda_j\,g_j^2+(\delta+\epsilon)\theta_s\,=\,q_\epsilon+
\delta\sum_{k=r+1}^s\sum_{j=1}^n\frac{x^2_i}{k{\rm
!}}\,=\,\hat{q}_\epsilon,\]
where $\hat{q}_\epsilon$ is a s.o.s. polynomial.

As $\delta$ was chosen to satisfy $\delta\leq \epsilon$, we obtain
\[f+\sum_{j=1}^m\lambda_j\,g_j^2+2\epsilon\theta_s\,=\,\hat{q}_\epsilon+
(\epsilon-\delta)\theta_s\,=\,\hat{\hat{q}}_\epsilon,\]
where again, $\hat{\hat{q}}_\epsilon$ is a s.o.s. polynomial. $\qed$

\section{Appendix}
In this section, we first briefly recall some basic results of linear
programming in infinite-dimensional spaces, and then 
present auxiliary results that are used in some of the proofs in \S\ref{proofs}.

\subsection{Linear programming in infinite dimensional spaces}
\label{LP}
\subsubsection{Dual pairs}
Let $\X,\Y$ be two arbitrary (real) vector spaces, and let $\ip{.}{.}$ be a bilinear form on 
$\X\times\Y$, that is, a real-valued function on $\X\times\Y$ such that

$\bullet$ the map $x\mapsto \ip{x}{y}$ is linear on $\X$ for every $y\in\Y$

$\bullet$ the map $y\mapsto \ip{x}{y}$ is linear on $\Y$ for every $x\in\X$.

Then the pair $(\X,\Y)$ is called a {\bf dual pair} if the bilinear form {\it separates} points in
$\X$ and $\Y$, that is,

$\bullet$ for each $0\neq x\in\X$, there is some $y\in\Y$ such that $\ip{x}{y}\neq 0$, and

$\bullet$ for each $0\neq y\in\Y$, there is some $x\in\X$ such that $\ip{x}{y}\neq 0$.

Given a dual pair $(\X,\Y)$, we denote by $\sigma(\X,\Y)$ the {\bf weak topology} on $\X$
(also referred to as the $\sigma$-topology on $\X$), namely the coarsest - or weakest - topology on $\X$, under which all the elements of $\Y$ are continuous when regarded as linear forms $\ip{.}{y}$ on $\X$.

Equivalently, the base of neighborhoods of the origin of the $\sigma$-topology is the family of all sets of the form
\[N(I,\epsilon)\,:=\,\{x\in\X\,\vert\quad \vert \ip{x}{y}\leq\epsilon,
\quad \forall y\in I\},\]
where $\epsilon>0$ and $I$ is a {\it finite} subset of $\Y$. (See for 
instance Robertson and Robertson \cite[p. 32]{robertson}.)
In this case, if $\{x_n\}$ is a net or a sequence in $\X$, then $x_n$ 
{\it converges} to $x$ in the weak topology $\sigma(\X,\Y)$ if
\[\ip{x_n}{y}\,\to\,\ip{x}{y},\qquad \forall y\in\Y.\]

\begin{definition}
\label{def1}
Let $(\X,\Y)$ and $(\Z,\W)$ be two dual pairs of vector spaces, and 
$G:\X \to \Z$, a linear map.

{\rm (a)} $G$ is said to be {\bf weakly continuous} if it is continuous with respect to the weak topologies $\sigma(\X,\Y)$ and $\sigma(\Z,\W)$; that is, if $\{x_n\}$ is a net in $\X$ such that $x_n\to x$ in the weak topology $\sigma(\X,\Y)$, then 
$Gx_n\to Gx$ in the weak topology $\sigma(\X,\Y)$, i.e.,
\[\ip{Gx_n}{v}\,\to\,\ip{Gx}{v},\qquad \forall v\in\W.\]

{\rm (b)} The {\bf adjoint} $G^*:\W\to\Y$ of $G$ is defined by the relation
\[\ip{Gx}{v}\,=\,\ip{x}{G^*v},\qquad \forall x\in\X,\,v\in\W.\]
\end{definition}
The following proposition gives a well-known (easy to use) criterion for the map
$G$ in Definition \ref{def1}, to be weakly continuous. 

\begin{prop}
\label{proplp}
The linear map $G$ is weakly continous if and only if its adjoint $G^*$ maps 
$\W$ into $\Y$, that is, $G^*(\W)\subset \Y$.
\end{prop}

\subsubsection{Positive and dual cones}
Let $(\X,\Y)$ be a dual pair of vector spaces, and $C$ a {\it convex cone}
in $\X$, that is, $x+x'$ and $\lambda x$ belong to $C$ whenever 
$x$ and $x'$ are in $C$ and $\lambda>0$. Unless explicitly stated otherwise, 
we shall assume that $C$ is not the whole space, that is, $C\neq\X$, and 
that the origin (the zero vector in $\X$) is in $C$. In this case, $C$ defines a partial order $\geq$ in $\X$, such that
\[x\,\geq x'\quad\Leftrightarrow\quad x-x'\,\in C ,\]
and $C$ is referred to as a {\it positive cone} in $\X$. The {\bf dual cone} of 
$C$ is the convex cone $C^*$ in $\Y$ defined by
\[C^*\,:=\,\{y\in\Y\,\vert\quad \ip{x}{y}\,\geq\,0,\quad \forall x\in C\}.\]

\subsubsection{Infinite linear programming (LP)}
An infinite linear program requires the following components:

$\bullet$ two dual pairs of vector spaces $(\X,\Y)$.

$\bullet$ a weakly continuous linear map $G: \X\to\Z$, with adjoint 
$G^*:\W\to\Y$.

$\bullet$ a positive cone $C$ in $\X$, with dual cone $C^*$ in $\Y$; and

$\bullet$ vectors $b\in\Z$ and $c\in\Y$.

\noindent
Then the {\bf primal} linear program is
\begin{equation}
\label{primallp}
\P:\quad\begin{array}{l}
{\rm minimize}\:\ip{x}{c}\\
\mbox{subject to: }Gx\,=\,b,\quad x\in C.\end{array}
\end{equation}
The corresponding {\bf dual} linear program is
\begin{equation}
\label{duallp}
\P^*:\quad\begin{array}{l}
{\rm maximize}\:\ip{b}{w}\\
\mbox{subject to: }c-G^*w\,\in\,C^*,\quad w\in \W.\end{array}
\end{equation}
An element of $x\in\X$ is called feasible for $\P$ if it satisfies 
(\ref{primallp}), and $\P$ is said to be consistent if it has a feasible solution. 
If $\P$ is consistent then its value is defined as 
\[\inf\P\,:=\,\inf \,\{\ip{x}{c}\,\vert\,\quad x\mbox{ is feasible for }\P\};\]
otherwise, $\inf\P=+\infty$.
The linear program $\P$ is solvable if there is some 
feasible solution $x^*\in\X$, that achieves the value $\inf\P$; then $x^*$ is an optimal solution of $\P$, andf one then writes $\inf\P=\min\P$.
The same definitions apply for the dual linear program $\P^*$.

The next result can be proved as in elementary (finite-dimensional) LP.

\begin{prop}[Weak duality]
\label{lpduality}
If $\P$ and $\P^*$ are both consistent, then their 
values are finite and satisfy $\sup\P^*\leq\inf\P$.
\end{prop}

There is no {\bf duality gap} if $\sup\P^*=\inf\P$, and strong duality holds
if $\max\P^*=\min\P$, i.e., if there is no duality gap, and both $\P^*$ 
and $\P$ are solvable.

\begin{theorem}
\label{nash}
Let $D$ be the set in $\Z\times\R$, defined as
\begin{equation}
\label{setd}
D\,:=\,\{(Gx,\,\ip{x}{c})\,\vert\quad x\in C\}.
\end{equation}
If $\P$ is consistent with finite value, and $D$ is weakly closed (i.e., closed in the weak topology 
$\sigma(\Z\times\R,\W\times\R)$), then $\P$ is solvable and there is no
duality gap, i.e., $\sup\P^*=\min\P$.
\end{theorem}
(See Anderson and Nash \cite[Theor. 3.10 and 3.22]{nash}.

\subsection{Auxiliary results}
\label{aux}
Let $\mathcal{B}$ be the Borel
sigma-algebra of $\R^n$, $\co$ be the space of bounded continuous
functions that vanish at infinity, and let $\theta_r$ be as in
(\ref{theta}). Let $M(\R^n)$ be the space of finite signed Borel measures on
$\R^n$.

\begin{lemma}
\label{lemmamom}
Let $r\geq1$, and let $\{\mu_j\}_{j\in J}\subset M(\R^n)$ be a sequence of 
probability measures, such that 
\begin{equation}
\label{lemmamom-1}\sup_{j\in J} \,\int \theta_r\,d\mu_j\,<\,\infty.
\end{equation}
Then there is a subsequence $\{j_k\}\subset J$ and a probability measure
$\mu$ on $\R^n$ (not necessarily in $\mathcal{M}$), such that
\[\lim_{k\to\infty}\,\int f\,d\mu_{j_k}\,=\, \int f\,d\mu,\]
for all continuous functions $f:\R^n\to \R$, such that $f/\theta_{r-1} \in\cb$.
\end{lemma}

\begin{proof}
$\theta_r$ is a {\it moment} function (see (\ref{momfunc})), and so,
(\ref{lemmamom-1}) implies that the
sequence $\{\mu_j\}$ is tight. Hence, as $\R^n$ is a metric space,
by Prohorov's Theorem \cite[Theor. 1.4.12]{hernand}, 
there is a subsequence $\{j_k\}\subset J$ and a measure
$\mu\in M(\R^n)$ such that $\mu_{j_k}\Rightarrow \mu$, i.e.,
\begin{equation}
\label{weak}
\int h\,\mu_{j_k}\,\to\, \int h\,d\mu ,
\end{equation}
for all $h\in\cb$, with $\cb$ being the space of {\it bounded continuous} functions $h:\R^n\to\R$.
Next, let $\nu_{j_k}$ be the measure obtained from
$\mu_{j_k}$ by:
\[\nu_{j_k}(B)\,:=\,\int_B \,\theta_{r-1}\,d\mu_{j_k},\qquad B\in\mathcal{B}.\]
Observe that from the definition of $\theta_r$, the function
$\theta_r/\theta_{r-1}$ is a moment function, for every $r\geq 1$.
And one has,
\[\sup_{k}\,\int \theta_r /\theta_{r-1}\,d\nu_{j_k}\,=\,
\sup_k\int \theta_r\,d\mu_{j_k}\,<\,\infty,\]
because of (\ref{lemmamom-1}). 
Observe that $\nu_{j_k}(\R^n)\leq \rho$ for all $k$, and so, we may
consider a subsequence of $\{j_k\}$ (still denoted $\{j_k\}$ for
simplicy of notation) such that $\nu_{j_k}(\R^n)\to \rho \,(>0)$ as $k\to\infty$.
With $\hat{\nu}_{j_k}:=\nu_{j_k}/\nu_{j_k}(\R^n)$, for all $k$, it
follows that the sequence of probability measures $\{\hat{\nu}_{j_k}\}_k$
is tight, which implies that there is a subsequence $\{j_n\}$ of
$\{j_k\}$, and a measure $\hat{\nu}\in M(\R^n)$,  such that
\[\mbox{ as }n\to\infty,\quad\int h\,d\hat{\nu}_{j_n}\,\to\, \int
h\,d\hat{\nu}, \quad 
\forall h\in\cb.\]
Since $\nu_{j_k}(\R^n)\to \rho$ as $k\to\infty$, we immediately get
\[\int h\,d\nu_{j_n}\,=\,\int h\,(\rho+\nu_{j_n}(\R^n)-\rho)\,d\hat{\nu}_{j_n}
\,\to\, \int h\rho \,d\hat{\nu}, \quad \mbox{ as }n\to\infty,\]
for all $h\in\cb$. Equivalently, with
$\nu:=\rho\hat{\nu}$,
\begin{equation}
\label{fin}
 \mbox{ as }n\to\infty,\quad\int h\,d\nu_{j_n}\,\to\, \int
h\,d\nu, \quad\forall h\in\cb.
\end{equation}
But as $h/\theta_{r-1}\in\cb$ whenever $h\in\cb$, (\ref{fin}) yields
\[\int h/\theta_{r-1}\,d\nu\,=\,
\lim_{n\to\infty}\int h/\theta_{r-1}\,d\nu_{j_n}\,=\,
\lim_{n\to\infty}\int h\,d\mu_{j_n}\\
\,=\,\int h\,d\mu,\]
for all $h\in\cb$. 

As both $\mu$ and $\theta_{r-1}^{-1}d\nu$ are finite measures, this implies that 
\begin{equation}
\label{defnu}
\mu(B)\,:=\,\int_B \,(1/\theta_{r-1})\,d\nu,\qquad B\in\mathcal{B}.
\end{equation}
As the subsequence $\{j_n\}$ was arbitrary, it thus follows that
the whole subsequence $\{\nu_{j_k}\}$ converges weakly to $\nu$.

Next, let $f:\R^n\to \R$ be continuous and such that
$f/\theta_{r-1}\in\cb$. As $k\to\infty$, from (\ref{fin}), 
\[\int (f/\theta_{r-1})\,d\nu_{j_k}\,\to\,
\int (f/\theta_{r-1})\,d\nu,\]
and so, 
\begin{eqnarray*}
\int \,f\,d\mu_{j_k}&=&\int (f/\theta_{r-1})\,\theta_{r-1}
\,d\mu_{j_k}\,=\,\int (f/\theta_{r-1})\,\,d\nu_{j_k}\\
&\to& \int (f/\theta_{r-1})\,d\nu\,=\,\int \,f\,d\mu,\quad \mbox{[by (\ref{defnu})],}
\end{eqnarray*}
the desired result.
\end{proof}
\begin{lemma}
\label{newlemma}
Let $\mu$ be a measure on $\R^n$ (with $\mu(\R^n)=y_0$) be such that
\begin{equation}
\label{newlemma-1}
\sup_{i=1,\ldots,n}\:\sup_{0\leq k\leq r}\:\int x^{2k}_i\,d\mu\,\leq\ S.
 \end{equation}
Then, 
\begin{equation}
\label{2}
\sup_{\alpha\in\N^n; \vert\alpha\vert \leq r}\:\vert\int
x^{\alpha}\,d\mu\,\vert\,\leq\ \sqrt{y_0S}.
 \end{equation}
\end{lemma}
\begin{proof}
Let $\y=\{y_\alpha\}_{\vert\alpha\vert\leq 2r}$, be the sequence
of moments, up to order $2r$, of the measure $\mu$, and let
$M_r(\y)$ be the moment matrix 
defined in (\ref{mommatrix}), associated with $\mu$. Then,
(\ref{newlemma-1}) means that those diagonal elements of $M_r(\y)$, denoted
$y^{(i)}_{2k}$ in Lasserre \cite{lasserre}, are all bounded 
by $S$. Therefore, by Lemma 6.2 in \cite{lasserre},
{\it all} diagonal elements of $M_r(\y)$ are also bounded by $S$, i.e.,
\begin{equation}
\label{conseq}
y_{2\alpha}\,\leq\, S,\quad\forall \alpha\in\N^n,\:\vert\alpha\vert\leq r,
\end{equation}
and so are all elements of $M_r(\y)$ (because $M_r(\y)\succeq 0$).
Next, consider the two columns (and rows) 
$1$ and $j$, associated with the
monomials $1$ and $x^\alpha$, respectively, and with
$\vert\alpha\vert\leq r$, that is, $M_r(\y)(1,1)=y_0$ and $M_r(\y)(1,j)=y_\alpha$.
As $M_r(\y)\succeq 0$, we immediately have
\[M_r(\y)(1,1)\times M_r(\y)(j,j)\,\geq\,M_r(\y)(1,j)M_r(\y)(j,1)\,=\,M_r(\y)(1,j)^2,\]
that is, $y_0y_{2\alpha}\geq y_\alpha^2$. Using that
$\vert\alpha\vert \leq r$ and (\ref{conseq}), we obtain
$y_0S\geq y_\alpha^2$, for all $\alpha, \:\vert\alpha\vert\leq r$, the
desired result (\ref{2}).
\end{proof}

\bibliographystyle{amsplain}

\end{document}